\documentclass{article}
\usepackage{graphicx}

\usepackage{amsmath,amssymb}
\usepackage{amsthm}
\usepackage{xcolor}
\usepackage[T1]{fontenc}
\usepackage[utf8]{inputenc}
\usepackage{indentfirst}

\usepackage{authblk}

\usepackage{hyperref}

\theoremstyle{plain}

\newtheorem*{theorem*}{Theorem}

\theoremstyle{definition}

\newtheorem*{problem*}{Problem}

\theoremstyle{remark}

\title{PRIMES STEP Experience}
\author{Slava Gerovitch}
\author{Tanya Khovanova}
\affil{MIT}
\date{}

\begin{document}

\maketitle

\begin{abstract}
PRIMES STEP is a mathematical outreach program established at MIT in 2015. STEP students study advanced topics beyond the school curriculum and conduct group research projects, often leading to publication. This article discusses the program's history, admissions process, lesson organization, interactive teaching style, and teamwork, and provides advice on how to choose research projects, encourage students, and keep them engaged. This paper would be useful for math teachers and instructors in after-school math programs.
\end{abstract}

\section{Introduction}

Who could have guessed that a group of 7th-9th graders would be able to find a connection between the Stable Matching problem and Sudoku? However, in 2021, this connection was discovered by participants in the PRIMES STEP program at MIT. The results came out in two papers: `The Stable Matching Problem and Sudoku', published in The College Mathematics Journal~\cite{SMPAndSudoku}, and `Sequences of the Stable Matching Problem', published in the Journal of Integer Sequences~\cite{SMPSequences}.

How was this breakthrough achieved by such young students? We will look into the history of STEP and share the experience of running this program, in the hope of stimulating math research activities in schools and after-school programs to inspire young mathematicians.

\section{Program structure}

\subsection{History}

\textbf{Slava:} Let us recall how it started. The idea of STEP emerged after our PRIMES program had become very successful. The MIT Department of Mathematics started PRIMES (Program for Research in Mathematics, Engineering, and Science for High School Students) in 2011 with 21 local students doing research in mathematics and computational biology \cite{EGK}. Our team included Pavel Etingof as the Chief Research Advisor, you as Head Mentor, and me as Program Director. The following year, we added computer science projects for local students and continued expanding the program every year, both in terms of the number of students and the range of research options. In 2013, we opened a nationwide remote math research program (PRIMES-USA) and an enrichment math program for local students (PRIMES Circle); in 2015, we launched a two-week mathematical talent accelerator residential summer program (MathROOTS), and in 2016, jointly with the Art of Problem Solving, we set up a massive collaborative online research forum open to all high school and college students around the world (CrowdMath). Currently, PRIMES serves more than 150 students annually and is unique among educational outreach programs in its rigor, duration, population it serves, and successes it has helped students achieve. Since the inception of the program, almost 450 papers have been posted on the PRIMES website; more than 90 of them have already been published in academic journals and professional conference proceedings. PRIMES students have won almost 400 awards in national and international science competitions.
E
\textbf{Tanya:} Pavel once jokingly suggested expanding PRIMES to kindergarten. We laughed, but then we started discussing whether math research could be done in middle school, and I was tasked to do it. I was afraid at first that it would be impossible. But then I started brainstorming. I made several decisions that made it work.

First, we decided that it should be group research. The students are too young to do individual research. I was also scared of what could happen if one student lost interest. If they are in a group, they can still pick one part of the project that is more appealing to them.

The second decision was that the program should not concentrate on research only; a big chunk of it should be a math club. The reason was that even at PRIMES, not every student wanted to be a mathematician. PRIMES STEP students are younger, so they do not know what they want to do when they grow up. It would be unwise to force them to do only research. When research is only a part of it, then they learn a lot of mathematics. And research is useful to teach them how to work on a project for a long time, plus it gives them a taste for research in general. Some students later confirmed that the program taught them more discipline and better work habits.

And third, I made a decision about the material to cover. The idea was not to overlap much with a school program. If I teach them something they know, they will be bored at STEP. If I teach them something that they later learn at school, they will be bored at school. Thus, my material is very diverse and not mainstream.

\subsection{Goal}

\textbf{Slava:} We decided to name the program PRIMES STEP (Solve–Theo\-rize–Explore–Prove), suggesting that this is a step into the world of mathematical discovery. We advertised it as a program for curious middle-schoolers who like solving challenging problems and are ready to explore advanced math topics. What do you expect your students to learn in the program? Is it a body of mathematical knowledge or a set of skills? What is STEP for?

\textbf{Tanya:} My goal in life and everywhere is whatever I do is to go for the biggest impact. In the case of STEP, it also means I not only teach mathematics, but I also teach them how to think and how to work. The goal of STEP is to use mathematics to expand the students' minds and to encourage their creative thinking. We do not repeat the middle-school math curriculum. Instead, we ask students to solve tricky problems, logical puzzles, and fun brain teasers, as well as Olympiad-style problems and questions that might be asked at a Wall Street interview. STEP students also emulate the research process by inventing new problems and being the first to solve them. We even learn an occasional magic trick. The main goal is to have fun!

\subsection{Outcome}

\textbf{Slava:} The first year, there was a group of 11 students in STEP, and they wrote two articles, both now published in \textit{Recreational Mathematics Magazine} \cite{STEPWhosGuilty} and \textit{Math Horizons} \cite{STEPAlternator}. In subsequent years, we had two groups of ten students each. 

\textbf{Tanya:} We call them senior and junior groups, with more advanced students in the former. Each group produces 1-2 papers per year. 

\textbf{Slava:} What is the range of topics covered by STEP papers?

\textbf{Tanya:} Most of our papers can be viewed as research in recreational mathematics. The topics range from logic to coin weighings to games to magic squares. Some of our papers moved into modern areas of research: combinatorics, number theory, information theory, game theory, and geometry. We post all our papers on arXiv.org and provide links on the STEP website.

\textbf{Slava:} Are these papers subsequently published?

\textbf{Tanya:} Most papers are either already published or submitted to journals. Due to the blind peer review procedure in some journals, reviewers are unaware that the authors are very young students, so we can be sure that the paper is approved for publication solely on merit.

\textbf{Slava:} What was the highest STEP achievement over the years? 

\textbf{Tanya:} One year was especially prolific when the senior group produced four papers. We started with the game of EvenQuads~\cite{LR}. The first paper~\cite{STEPCardGames} describes the connections of EvenQuads to other card games: SET, Socks, and Spot It!. We studied analogs of magic squares for EvenQuates and wrote about them in~\cite{STEPQuadSquares}. In~\cite{STEPMaxNumQuads}, we studied the maximum possible number of quads in a set of $\ell$ cards depending on the deck size. Our most advanced paper connects the game of EvenQuads to error-correcting codes and is published in Springer Nature Computer Science.

\subsection{Admission}

\textbf{Slava:} What kind of students are you looking for?

\textbf{Tanya:} We are looking for students who are interested in mathematics and want to work on research. But realistically, we sometimes have students who do not know that our program conducts research. They apply because some of their friends recommended it, or their parents want them to be here for their future resumes. 

\textbf{Slava:} Alright! How do you select students for STEP? 

\textbf{Tanya:} The only selection method we have is an entrance test. We do need students who are strong in mathematics. After the entrance test, I often interview the students. This is my chance to tell them that they will do research and see their reaction. Once, I interviewed a student who seemed not interested. His parents brought him to the entrance test and never explained what the program was. He thought that this was another RSM (Russian School of Mathematics). After we talked, he got excited and grew into a great student.

\textbf{Slava:} And what kind of student succeeds in STEP?

\textbf{Tanya:} Most of them succeed. They enjoy the material. They help each other. They become friends. They want the paper to be out there, so they're excited about the class.

Very rarely do I have students who procrastinate. Many years ago, I had a student who didn't contribute anything to research. I told him that I would remove his name from the paper unless he did something. He helped prepare beautiful slides for the conference, so I kept his name.

Now, I give research points, and I have a cutoff for a person to become a coauthor. Since then, I have had only one student who did not submit research and did not do much of the homework. I asked this student to leave the program after the first semester.

\section{Class structure}

\textbf{Slava:} Let us get to the nitty-gritty details of how the program works. What is a STEP session like?

\textbf{Tanya:} The class runs for two hours with a break. It is usually divided into four parts: 
\begin{enumerate}
    \item homework discussion,
    \item hands-on activity,
    \item new mathematical topic, and
    \item research discussion.
\end{enumerate}

\subsection{Homework discussion}

\textbf{Slava:} Let's start with the homework assignment. Do you assign homework? Is it like regular school homework?

\textbf{Tanya:} My homework has about 10 questions. It starts with a warm-up. As I teach the students how to think, I try to be broader than mathematics. Thus, in a warm-up, I might include a riddle, a joke, a physics puzzle, a logic puzzle, a linguistic puzzle, a tricky question, and so on.

Here is an example of a physics problem followed by a riddle.

\begin{problem*}
Where on Earth would an object weigh the smallest? The largest?
\end{problem*}

\begin{problem*}
What kind of frog can jump higher than a building?
\end{problem*}

It could even be a joke question. One homework problem, for example, was the following.

\begin{problem*}
Do not read this sentence.
\end{problem*}

My favorite answer was: Do not read this answer.

The core of the homework is questions related to the current topic. Sometimes, I also have questions from past math competitions for training. Sometimes, I have random math questions for fun. I end my homework with one or two challenge questions. For example, when we studied the Catalan numbers, the homework contained the following challenge question.

\begin{problem*}
Find a formula for the number of Dyck words of length $2n$ starting with two Xs and ending with two Ys.
\end{problem*}

When we studied how integers are represented in computers, I gave them the following challenge.

\begin{problem*}
In your favorite programming language, write a program that, when run, will print out its own source code. In addition, consider English as a programming language: write a command in English to output the text of the command.
\end{problem*}

During the discussion of this homework, I explained to them that they wrote a quine. Nowadays, I should ask the students to write a ChatGPT request that outputs their request back exactly as it was.

Occasionally, I give them a special homework: a puzzle from the MIT Mystery Hunt. Puzzles from puzzle hunts are different from regular problems. Usually, the solvers are not told what to do; they need to figure it out themselves. I give one puzzle for everyone to solve together. It is great training for pattern recognition. Also, it builds teamwork. For example, for my class on cellular automata, I gave them the puzzle The Next Generation, which can be found at \url{https://puzzles.mit.edu/2018/full/puzzle/the_next_generation.html}. 

So, the request is that the homework be individual work because the homework is for practice. I grade the homework. I explain to them that the value of my grading is rather feedback than anything else. I do not use the grades anywhere except maybe in recommendation letters.

\textbf{Slava:} How do you discuss homework in class? Do you go problem by problem?

\textbf{Tanya:} I usually start my class by discussing the homework. I do not discuss every problem. If the solution is straightforward and everyone has solved it, there is no reason to discuss it. If one student didn't solve it, and the solution is available online, I can skip it too. I only discuss interesting ideas they suggested or missed. We also talk about the big picture.

Here is a homework problem borrowed from Alexander Karabegov.

\begin{problem*}
A function $f: \mathbb{R} \to \mathbb{R}$ satisfies the following conditions:
\begin{itemize}
\item $f(1) = 1$.
\item $f(-x) = -f (x)$ for all $x$; that is, $f(x)$ is an odd function.
\item $f(1 - x) = f(1 + x)$ for all $x$; that is, the function $g(x) = f(1 -x)$ is even.
\end{itemize}
Calculate $f(2024)$ and $f(2025)$.
\end{problem*}

Everyone solved it, but only one student noticed that the problem was inspired by a sine function. So, we discussed that.

After we finish the homework, I often present a toy or a trick.

\subsection{Hands-on activity}

\textbf{Slava:} Do students work only with pen and paper?

\textbf{Tanya:} No, not all! I think for younger students, it is better if they can see, play, and touch physical objects. I bring a lot of such objects to class. I often use a toy or a magic trick for a short transition from a homework discussion to a new topic.

\textbf{Slava:} You wrote on your blog that you hate crocheting. However, you crocheted a lot of mathematical objects. Do you use them in class?

\textbf{Tanya:} I wanted to teach hyperbolic surfaces, so I crocheted ten of them. It is much easier to understand them if you can hold a surface, fold it, find lines on it, and so on. I tried to buy them, but couldn't find them anywhere, so I made them. When the students hold a crocheted hyperbolic surface, it is easy to show how geodesics work. You can flatten out a small neighborhood of a point and see how close it is to a plane piece. Sometimes, one can find a line between two points on a crocheted surface by picking the points and stretching the fabric. If the given surface is a hyperbolic surface of constant curvature, you can find lines by folding. The crease corresponds to a geodesic. It is fun to prove that the parallel postulate is independent of other axioms by fiddling with a physical object.

\textbf{Slava:} Do physical objects help mathematicians think?

\textbf{Tanya:} Mathematics is about ideas; however, people have a variety of senses to help them remember, understand, and process ideas. A lot of things are easier to understand while playing with physical objects. For example, it is much easier to calculate the linking number of a link by fiddling with the link rather than looking at its picture \cite{TK}. 

\textbf{Slava:} What other objects do you bring to class?

\textbf{Tanya:} I am the proud owner of a dozen Rubik's cubes, enough for all students. When I teach group theory, I bring the cubes to the class and explain that the moves of the cube form a group. I explain to the students that a commutator of two elements is a measure of non-commutativity in a group. One might hope that a commutator of two group elements is a group element that is not far from the identity. If we take a commutator of two different 90-degree rotations of the neighboring faces of the Rubik's cube, the result shuffles four small vertex cubes and three small edge cubes. This allows me to explain one of the ways to solve the cube. It is a long way, but it helps with explaining groups.

\textbf{Slava:} I can see that math can be fun, and fun can be math. Any other ways to combine the two?

\textbf{Tanya:} Once, I brought two inflatable globes to class. I was teaching a theorem that any movement in 3D space that keeps the origin and orientation is a rotation. Equivalently, any sphere rotation has to have two antipodal fixed points. So, we played a game. I fixed one globe and asked the students to rotate the other globe randomly and then find points on Earth in the same relative position. I found this game in the book \cite{PM21}.

\textbf{Slava:} What was the funniest game you played?

\textbf{Tanya:} I like reenacting famous puzzles. Consider the following hat puzzle invented by Konstantin Knop and Alexander Shapovalov. It appeared (in a different wording) in March 2013 at the Tournament of the Towns:

\begin{quote}
    A sultan decides to give 100 of his sages a test. The sages will stand in line, one behind the other, so that the last person in the line sees everyone else. The sultan has 101 hats, each of a different color, and the sages know all the colors. The sultan puts all but one of the hats on the sages. The sages can only see the colors of the hats on people in front of them. Then, in any order they want, each sage guesses the color of the hat on his own head. Each hears all previously made guesses, but other than that, the sages cannot speak. They are not allowed to repeat a color that was already announced. Each person who guesses his color wrong will get his head chopped off. The ones who guess correctly go free. The rules of the test are given to them one day before the test, at which point they have a chance to agree on a strategy that will minimize the number of people who guess their hat color incorrectly. What should that strategy be?
\end{quote}

There is a good strategy, using modulo arithmetic, where no more than three people fail to eventually guess their hat color. But the best strategy, where 99 people are guaranteed survival, uses permutations. The students love reenacting hat puzzles. This forces students to work as a team, and they enjoy explaining the guessing algorithm to their friends.

\textbf{Slava:} Are you playing math games just for fun, or are they related to research?

\textbf{Tanya:} Of course, games are always related to research. Once, we started the academic year by learning to play several card games, including SET, Socks, EvenQuads, and Spot It! We formalized the games in the language of linear algebra and found how they are connected. In particular, we discovered a way to use the cards from one game to play other games. In the end, we wrote a lovely research paper comparing the four games \cite{STEPCardGames}.

\textbf{Slava:} You have mentioned that you use magic tricks. Is there room for magic in mathematics?

\textbf{Tanya:} There are tons of magic tricks that are not sleight of hand but are based on math. Tricks are fun, and I use them whenever I can.

\textbf{Slava:} Can you give an example?

\textbf{Tanya:} Sure. I ask each student to write down a digit that is not zero or one. I tell them that they will learn later why zeros and ones should be excluded. I ask them to multiply all their digits without showing them to me. Next, they must remove the first digit of the result and write the rest of the digits of the result on the blackboard in random order. Then, I say ``abracadabra'' and guess the first digit.

I use this trick when we study divisibility or probability. The point is that, with a very high probability, the product of the digits is divisible by 9, which allows me to quickly figure out the missing digit.

The only sleight of hand here is my explanation of why the ones and zeros should be excluded. Zeros should obviously be excluded, as the product would become zero. The ones are boring, I say, since when you multiply the digits, they do not affect the product. However, the real reason for excluding the ones is to increase the probability of divisibility by 9.

My students enjoy figuring out this trick. They get excited and start calculating the probability that it works. By the way, what is this probability if there are ten students in a group?

\subsection{New mathematical topic}

\textbf{Slava:} How do you choose topics for study?

\textbf{Tanya:} At the beginning of the academic year, I often teach new topics that might be useful in research. For example, the Fall of 2023 started with Fibonacci and Lucas numbers. We continued with some number theory and sequences for several classes because we needed this for research. However, I wanted to present various topics, and as soon as we finished with sequences, we moved to logic gates and had several classes related to computer science. In the next semester, I got excited about the recent development in aperiodic tiling, and we had a class on that.

The junior project needed some discussion on information theory. We studied information theory using coin and hat puzzles. Of course, I brought hats to class, and we reenacted the puzzles. The topics could be fair division, dancing tangles, Fermat's last theorem, or something else. 

\textbf{Slava:} Fermat's last theorem is very famous and took hundreds of years to prove, right? How did your class go?

\textbf{Tanya:} I wrote the theorem's statement and the word ``Proof.'' 

\begin{theorem*}[Fermat's Last Theorem.]
No three positive integers $(a, b, c)$ can satisfy the equation $a^n + b^n = c^n$ for any integer value of $n$ greater than 2. 
\end{theorem*}

\begin{proof}
\end{proof}

The senior group laughed at this point. 

My classes are very interactive, so we usually prove things together. I often pause and ask them what is next. So, we proved that Fermat's last theorem can be reduced to cases of $n$ being an odd prime and $n=4$. Then we reduced it to the existence of primitive triples (triples of numbers $(a, b, c)$ such that they are pairwise coprime). Then, we started Case 1 of $n=4$. The previous class was on Pythagorean triples, so we were able to prove this case. Then, we proceeded with Case 2 of $n$ being an odd prime. The challenge for me was to use the blackboard in such a way so that I would run out of available margins at the right time.

\textbf{Slava:} Do you reuse your material?

\textbf{Tanya:} Because my students continue to show up in subsequent years, I do not repeat the same topics. I can recycle topics, but it takes several years to do so. Also, I discovered that I am easily bored, so I enjoy designing new material. For example, recently, I taught Eulerian numbers, which I had never taught before, and, telling you a secret, I didn't know about them myself until I stumbled upon them in my own PRIMES project.

At the end of the class, we discuss research.

\subsection{Research discussion}

\textbf{Slava:} How do you choose a research topic?

\textbf{Tanya:} I always promise a research project for the senior group. I usually prepare a general research topic and allow the students to brainstorm particular directions. My motivation is as follows. An important part of being a researcher is the ability to ask questions. So, I need to train them to do so. The benefit is that they are more invested when they suggest a particular question themselves, which helps the project.

\textbf{Slava:} Can you give an example of a research project?

\textbf{Tanya:} In our first research, we introduced a new type of coin: the alternator. The alternator can pretend to be either a real or a fake coin (which is lighter than a real one). Each time it is put on a balance scale, it switches between pretending to be either a real coin or a fake one. We wrote a paper \cite{STEPAlternator} where we solved the following problem: You are given $N$ coins that look identical, but one of them is the alternator. All real coins weigh the same. You have a balance scale that you can use to find the alternator. What is the smallest number of weighings that guarantees that you will find the alternator? While working on this research, we discovered a beautiful sequence my students had never heard of: Jacobsthal numbers.

\textbf{Slava:} And for the junior group?

\textbf{Tanya:} I always say that I do not guarantee a research project for the junior group, but we always end up having one. I always include games and magic tricks in the first few lessons and ask the students to invent questions. In one year, it took two months for the students to come up with a question, but it usually goes much faster.

\textbf{Slava:} What is the benefit of doing research in middle school?

\textbf{Tanya:} Many of our students shine at their schools. Everything is easy for them. I remember what my son once told me. He told me that school was easy for him as a breeze; the disadvantage of being fast and gifted is that you do not need to learn to work. One of the things our program achieves is teaching the kids how to work.

\textbf{Slava:} How is the work on a project organized?

\textbf{Tanya:} We discuss research in class, and I also send it as homework. At the beginning of the year, we concentrate on the background of the topic and potential ideas for research. Then, we start implementing the ideas and discussing theorems and proofs. In the end, we decide which results fit together and write a paper. Editing a paper is also a lot of work. 

\textbf{Slava:} Formal mathematical writing might seem boring to young students. How do you deal with that?

\textbf{Tanya:} Our papers have different flavors. If our research is suitable for a research journal, we use formal mathematical writing. It might seem boring, but it teaches precision and is very beneficial for students. Sometimes, we write recreational papers, and the writing style is different. For example, our paper \textit{The Struggles of Chessland} \cite{Chessland} was written as a fairy tale using a dialogue between a spoiled chess queen and her minions.

\textbf{Slava:} How do students work together on a joint paper?

\textbf{Tanya:} I usually send them a file with the current material and the questions highlighted. They are expected to pick several questions they like and answer them. I do not want all the students to answer the same question, as it is wasteful, but having 2 or 3 students answer the same question is useful. As they are prone to mistakes, the same answer from several students makes it more reliable. I also like that they can pick a question closer to their hearts. They can use questions according to their preferences: making calculations, writing proofs, developing new ideas or concepts, programming, and so on. The file slowly grows into a paper.

\textbf{Slava:} How do you make sure that they do research themselves rather than watching you do research?

\textbf{Tanya:} For most of the year, I try not to think about our research project myself. Every week, the students submit their research ideas, and I review them. Many students do not know how to write mathematically. They invent their notation on the spot and do not explain it. They skip steps in their thought process. I have to go through their submissions and find ideas that might be useful. Reviewing their research homework is the most time-consuming part of my work. The good news is that with time and my feedback, their writing improves, and the review becomes easier.

\textbf{Slava:} How does a research paper emerge from this process?

\textbf{Tanya:} I compile research submissions into a PDF file of our draft paper. I reformat the students' ideas into statements, proofs, examples, and remarks. I check what is missing and add my questions to the file. My questions are highlighted and are easy to find. I ask them to provide proofs for new statements, clarify steps in the existing proofs, calculate examples, suggest new ideas, check literature, brainstorm more potential research questions, and so on.

\textbf{Slava:} Is everything included in the paper?

\textbf{Tanya:} We often have tons of preliminary calculations that are not needed in the final version. For example, a student calculated the time complexity of one of the algorithms. It was a good exercise, but our paper was not about algorithms, so we did not include it.

\textbf{Slava:} Then the students present their work at the PRIMES spring conference in May. How do they prepare their talk?

\textbf{Tanya:} While I participate in paper writing, I do not participate in the preparations for their talk; I just provide feedback during a practice run. It is their presentation, and they are very proud of it. They go the extra mile for their talk. They often present their project as an engaging story, adding puns, magic tricks, and jokes. Some groups even turned it into a theater-like performance. Watching them is a lot of fun!

\section{Advice, issues, and pitfalls}

\subsection{Encouragement}

\textbf{Slava:} How do you maintain students' interest throughout the year?

\textbf{Tanya:} Most of them enjoy homework and try to solve all the problems. My homework is difficult, so not every class I have a student who solves everything. I announce the top scorers, and they are very proud of it. I also award points for research. One must accumulate a certain number of research points during the year to become a coauthor of the paper. If someone falls behind, I would talk to them, provide encouragement, figure out their interests, and tailor specific questions to their interests.

Also, at the end of each semester, I give out small math games and puzzles as prizes for the students who obtained many research points. During class, I award gold stars for good ideas. My current method came about during Covid, when we held sessions over Zoom. The first time a student online had a great idea, I took a gold star sticker and tried giving it to them as usual, but hit the screen. Everyone started laughing. So I just glued the gold star to my forehead. Everyone loved this method so much that I keep doing it now in person. After class, I proudly walk to my parking garage with my face covered in stickers.

\textbf{Slava:} Different students may have different strengths. How do you make them work together on a team project?

\textbf{Tanya:} I have a lot of questions in my research homework. I ask them to pick a couple of questions and resolve them. That way, they pick questions according to their taste. For example, students who like computer programming pick questions that require coding, and so on. I like it when two or more students answer the same question. If they obtain the same result, it is more likely to be correct.

\textbf{Slava:} What if different students contribute unequally---some more, some less. How do you handle such situations?

\textbf{Tanya:} In all papers with multiple authors, contributions are not distributed evenly. The students with a lot of points enjoy research very much. Once, I had a student with a suboptimal number of points. I talked to him and discovered that he enjoys the programming part more than the mathematical part. So I started giving him specific programming tasks, and he easily gained enough points to be a coauthor and started enjoying his contributions.

\subsection{Interactive style of teaching}

\textbf{Slava:} How do you engage students in the learning process?

\textbf{Tanya:} I think it is very important to be interactive. When I teach, I always give students a chance to solve a problem themselves or prove a theorem together. Sometimes, I even suggest that they finish my definition. I often ask what they think I will discuss next. This helps develop the mathematical taste and the vision of the big picture.

\textbf{Slava:} Are there any special tricks?

\textbf{Tanya:} I have a trick I learned from the famous mathematician John Conway. I saw John give the same lecture several times, and every time, he would stare at the board at a particular moment, waiting for the audience to ``help him.'' Unfortunately, when he got older, the audience did not react, as it looked more and more like a senior moment. I adapted the trick: sometimes, when I teach a new topic, I pretend that I forgot what comes next. I look at the board in the middle of my presentation and wait for them to help me. However, I explicitly tell my students that my pauses are staged.

Sometimes I make an error on purpose to see if they are paying attention. Of course, I can make an inadvertent error as well. I am glad to report that they do correct me.

\textbf{Slava:} And what do they learn by doing this?

\textbf{Tanya:} When they participate in the discussion, they remember it better. Plus, they cannot just sit back and relax. This keeps them engaged when they know that we are making progress together.

\subsection{Group work}

\textbf{Slava:} Do you encourage or discourage a competitive spirit?

\textbf{Tanya:} I teach them to be less competitive. When they first join the program, if they know the answer to my question, they often shout it out to show how smart they are. I explain to them that the program's goal is to teach students how to think. If someone shouts the answer to the question immediately, other students do not have time to think. The rule is that if you know the answer, you need to give other students a chance to think about the question. They get this idea very fast and start raising their hands.

\textbf{Slava:} Are the students working as one group or divided into subgroups?

\textbf{Tanya:} Sometimes, when we discuss research or I give them a complicated problem, I ask them to work in smaller groups. I noticed that sometimes it's better to mix groups up, so I always have a deck of cards with me. (I have a deck of cards anyway, in case I need to show a trick.) For example, I can use the cards to create random groups. I pick five red and five black cards and shuffle them. They pick a card, and thus, I divide the ten students into two random groups of five. 

They even told me that if I assign them to random groups, they get to know each other better and make more friends.

\textbf{Slava:} I think someone has suggested once that aides to Republican and Democratic legislators on the Hill should be randomly assigned seats in a cafeteria so that they could mingle with one another, and that would help build a better understanding between different political factions.

\subsection{Competitions}

\textbf{Slava:} You have a lot of experience with math competitions. You were the second young woman ever to win gold at the International Math Olympiad. Are you using techniques for competition training in STEP?

\textbf{Tanya:} We sometimes have homework from past competitions just for practice. The mathematics we discuss in class might be useful in competitions. However, the program's goal is research, and its benefits for competitions are indirect. Though I participated in a lot of competitions myself, I stopped liking them.

\textbf{Slava:} Why did it happen?

\textbf{Tanya:} Mostly because the problems changed. When I was young, there was no Internet, which meant that people would reuse problem ideas in different competitions in different countries. As a result, the percentage of problems with interesting ideas was high. And now, if an interesting problem appears in a competition, an hour later, its solution is posted on the Internet. Thus, you cannot ever reuse the problems. And I think because of that, it is difficult to design fresh problems with new, interesting, elegant ideas.

Plus, my students are mostly from middle school. Their main competition is MATHCOUNTS, which is based on speed, and speed is useful. But calculations are more and more outsourced to computers, and training for speed is less and less valuable over time.

\subsection{Comparison with other programs}

\textbf{Slava:} Some advanced summer math camps also offer research experience to their participants. How is STEP different?

\textbf{Tanya:} Many summer camps and programs offer research experience. MIT offers to its undergraduates a class titled `Project Laboratory in Mathematics'. In this class, the students are given a project, which is an open-ended question, where they need to discover some patterns and write up their results in a paper. However, the questions have known answers. The students have research experience by reproducing research that has been done before. Our goal is different. Our goal is not just to have research experience, but to actually produce research no one has done before. In the end, we post our papers on the arXiv and send them for publication in a journal.

\textbf{Slava:} Does STEP focus on a single math topic or cover a wide range?

\textbf{Tanya:} I never wanted to have a program on a particular subject, like combinatorics or number theory. There are several summer camps that are devoted to number theory, and it's quite cool, but in STEP, by my estimate, not more than one out of ten people might later become a mathematician. Everyone else would choose another profession. Why do they need advanced number theory? I have several thematic modules during the term. For example, I can have a module of four classes on combinatorics, then a class on invariants, then a class on meta-solving, and then continue with a module on number theory mixed with cryptography. I decided that I would try to be as diverse as possible and cover many topics. I can have fair division and voting power, followed by triangle centers and scaling arguments in nature. I try to cover a lot of territory. For a research project, of course, we go much deeper.

\subsection{Teamwork}

\textbf{Slava:} You have a bunch of excited students in the room. How do you make them work together?

\textbf{Tanya:} As I mentioned, my teaching is interactive. So, I start a topic and begin asking questions. For example, we notice that the order of an element of a group divides the order of the group, and I ask them why. They would start discussing it and then work out a solution together. Sometimes, I point out how one student's idea was a development of another student's idea. I really enjoy it when they bounce ideas off each other, and, with time, they start to enjoy it too.

In our classroom, we have walls completely covered with blackboards. When I ask them to work on a problem, they just stand up and go to the blackboards, self-organizing into smaller groups.

\subsection{Choosing research projects}

\textbf{Slava:} Where do you find topics for research?

\textbf{Tanya:} This is difficult. So far, miraculously, I have managed to do it every year. I look at my past projects to see if there are any unsolved problems. I go around and talk to various people. Sometimes I pick up a project at the Gathering for Gardner or the MOVIES conferences.

\textbf{Slava:} Do students participate in the choice of topic?

\textbf{Tanya:} The difference between research and problem-solving is that in research, you choose your own questions. It is a good idea to train the students to ask questions. However, I do not want them to choose a topic completely by themselves.

\textbf{Slava:} Why not?

\textbf{Tanya:} For example, once a student suggested doing research on negative number bases. Negative bases were invented more than 100 years ago; thus, a lot has been done. We can't just blindly go into negative number bases; we need a more specific question. 

\textbf{Slava:} So, do you always pick a topic yourself?

\textbf{Tanya:} For the senior group, I prepare a topic in advance, where I feel there is potential for progress and publication in a research journal. Then, I ask the students to brainstorm possible questions. Then, we check the literature and discuss which questions seem more doable and interesting. After we choose a couple of directions, we start working, but we continue discussing possible new directions for expanding the project. I close the option of suggesting new questions somewhere in March, close to the end of the academic year.

\textbf{Slava:} And how is a topic chosen for the junior group?

\textbf{Tanya:} The policy in the junior group is a bit different. The junior group is usually less advanced, so my goal is not to publish in a research journal but rather in a recreational journal. Again, it is difficult to find a good question, so at the beginning of the academic year, I include some puzzles in the homework, where I see potential for branching out. For example, during the 2023--2024 academic year, I showed them two magic tricks related to each other, and they just took off. The junior group wrote a fun paper~\cite{STEPCardTricksAndInformation} connecting the card tricks to information theory. They invented a new trick on the way!

\textbf{Slava:} Do projects for the senior group require more background?

\textbf{Tanya:} Yes, and I cover the basic background as part of my class topics. When the background becomes too complicated for both senior and junior groups, I discuss the advanced details with the senior group only. For example, in the 2024--25 academic year, I prepared a research project on chip-firing. We started the academic year by studying the dollar game and chip-firing in our classes for both groups. I also used chip-firing to teach some other things in mathematics, such as equivalence classes, graph theory, sequences, and other things. We discussed more advanced questions with the senior group only.

\textbf{Slava:} Do you also discuss the background for the junior group project in both groups?

\textbf{Tanya:} Yes, whenever there is something in the research background for one group that might be fun for both groups, I make it a topic for a class. For example, in the 2023--2024 academic year, during the second semester, it became clear that the junior group project required some background in information theory, so I taught it to both groups.

\subsection{Toy research or real research?}

\textbf{Slava:} Many of your research topics fall under the rubric of ``recreational mathematics.'' Is this something that could be called ``toy research''? Or is it real mathematical research? How useful is this kind of research for the future career of a mathematician?

\textbf{Tanya:} What STEP students do \textbf{is} research. We are doing something no one has done before. When we investigate the question, the very process is a research process. Recreational mathematics means that you do not need too much advanced knowledge to understand it, yet it is possible to discover some new phenomena or construct an object that might be useful for real mathematics. For example, while working on a project, we often stumble upon new integer sequences, and then we submit those sequences to the Encyclopedia of Integer Sequences \cite{OEIS}. It's not exactly a breakthrough, but it is useful for the mathematical community to have such interesting sequences.

In 2023, David Smith, who was not a professional mathematician, invented the monotile, which people had been looking for for years. This is related to the aperiodic tiling of the plane, and many mathematicians have been trying to find the monotile in vain. It is quite possible that we can also stumble on something that might be useful in the future.

Many of our papers have been published, not just in recreational mathematics journals but in regular research journals, too. It \textbf{is} real research!

\subsection{The impact of STEP}

\textbf{Slava:} What is the long-term impact of the program on STEP students?

\textbf{Tanya:} My students fall in love with mathematics and continue with MIT PRIMES and other programs. They enjoy research so much that after completing our program, they begin seeking new research opportunities.

From time to time, I run into my former students on the MIT campus. It is always great to see them grow and mature as scholars. They tell me that STEP research was a great experience for them, that they feel they have contributed to the world of knowledge.

\textbf{Slava:} Do you often hear from the alumni of the program?

\textbf{Tanya:} I often receive postcards with words of gratitude. MIT Admissions also regularly sends me letters saying that such and such a student they accepted named me as the person who made the greatest impact on their lives.

\textbf{Slava:} That's so wonderful to hear. I hope you continue running this unique program and raising new generations of budding mathematicians!

\section{Acknowledgments}

We wish to thank Pavel Etingof for suggesting the idea of PRIMES STEP and for his continuing guidance and support. We are also grateful to the MIT Mathematics Department for supporting the program.

\end{document}